\documentclass[english,12pt]{smfart}
\usepackage{amssymb}
\usepackage{amsfonts}
\usepackage{amscd}
\usepackage[all]{xypic}

\swapnumbers

\def\limproj{\mathop{\oalign{lim\cr
\hidewidth$\longleftarrow$\hidewidth\cr}}}

\def\Var{{\rm Var}}

\let\got\mathfrak
\def\gP{{\got P}}

\def\AA{{\mathbf A}}

\def\CC{{\mathbf C}}

\def\FF{{\mathbf F}}

\def\LL{{\mathbf L}}

\def\NN{{\mathbf N}}

\def\PP{{\mathbf P}}
\def\QQ{{\mathbf Q}}

\def\ZZ{{\mathbf Z}}

\def\cL{{\mathcal L}}
\def\cM{{\mathcal M}}

\def\cO{{\mathcal O}}

\mathchardef\alphag="7C0B
\mathchardef\betag="7C0C
\mathchardef\gammag="7C0D
\mathchardef\deltag="7C0E
\mathchardef\varepsilong="7C22
\mathchardef\varphig="7C27
\mathchardef\psig="7C20
\mathchardef\zetag="7C10
\mathchardef\epsilong="7C0F
\mathchardef\rhog="7C1A
\mathchardef\taug="7C1C
\mathchardef\upsilong="7C1D
\mathchardef\iotag="7C13
\mathchardef\thetag="7C12
\mathchardef\pig="7C19
\mathchardef\sigmag="7C1B
\mathchardef\etag="7C11
\mathchardef\omegag="7C21
\mathchardef\kappag="7C14
\mathchardef\lambdag="7C15
\mathchardef\mug="7C16
\mathchardef\xig="7C18
\mathchardef\chig="7C1F
\mathchardef\nug="7C17
\mathchardef\varthetag="7C23
\mathchardef\varpig="7C24
\mathchardef\varrhog="7C25
\mathchardef\varsigmag="7C26
\mathchardef\Omegag="7C0A
\mathchardef\Thetag="7C02
\mathchardef\Sigmag="7C06
\mathchardef\Deltag="7C01
\mathchardef\Phig="7C08
\mathchardef\Gammag="7C00
\mathchardef\Psig="7C09
\mathchardef\Lambdag="7C03
\mathchardef\Xig="7C04
\mathchardef\Pig="7C05
\mathchardef\Upsilong="7C07

\newtheorem{theorem}[subsubsection]{Theorem}
\newtheorem{lem}[subsubsection]{Lemma}
\newtheorem{klem}[subsubsection]{Key-Lemma}
\newtheorem{cor}[subsubsection]{Corollary}
\newtheorem{prop}[subsubsection]{Proposition}

\newtheorem{sconjecture}[subsubsection]{Rationality conjecture (strong form)}
\newtheorem{wconjecture}[subsubsection]{Rationality conjecture (weak form)}

\theoremstyle{definition}

\newtheorem{example}[subsubsection]{Example}

\newtheorem{def-prop}[subsubsection]{Proposition-Definition}
\newtheorem{def-theorem}[subsubsection]{Theorem-Definition}

\theoremstyle{remark}
\newtheorem{remark}[subsubsection]{Remark}
\newtheorem{remarks}[subsubsection]{Remarks}

\theoremstyle{plain}

\numberwithin{equation}{subsection}

\def\boxit#1#2{\setbox1=\hbox{\kern#1{#2}\kern#1}%
\dimen1=\ht1 \advance\dimen1 by #1
\dimen2=\dp1 \advance\dimen2 by #1
\setbox1=\hbox{\vrule height\dimen1 depth\dimen2\box1\vrule}%
\setbox1=\vbox{\hrule\box1\hrule}%
\advance\dimen1 by .4pt \ht1=\dimen1
\advance\dimen2 by .4pt \dp1=\dimen2 \box1\relax}

\renewcommand{\theequation}{\thesubsection.\arabic{equation}}

\let\got\mathfrak
\def\AA{{\mathbf A}}

\def\CC{{\mathbf C}}

\def\FF{{\mathbf F}}

\def\LL{{\mathbf L}}

\def\NN{{\mathbf N}}

\def\PP{{\mathbf P}}
\def\QQ{{\mathbf Q}}

\def\ZZ{{\mathbf Z}}

\def\cL{{\mathcal L}}
\def\cM{{\mathcal M}}

\def\cO{{\mathcal O}}

\mathchardef\alphag="7C0B
\mathchardef\betag="7C0C
\mathchardef\gammag="7C0D
\mathchardef\deltag="7C0E
\mathchardef\varepsilong="7C22
\mathchardef\varphig="7C27
\mathchardef\psig="7C20
\mathchardef\zetag="7C10
\mathchardef\epsilong="7C0F
\mathchardef\rhog="7C1A
\mathchardef\taug="7C1C
\mathchardef\upsilong="7C1D
\mathchardef\iotag="7C13
\mathchardef\thetag="7C12
\mathchardef\pig="7C19
\mathchardef\sigmag="7C1B
\mathchardef\etag="7C11
\mathchardef\omegag="7C21
\mathchardef\kappag="7C14
\mathchardef\lambdag="7C15
\mathchardef\mug="7C16
\mathchardef\xig="7C18
\mathchardef\chig="7C1F
\mathchardef\nug="7C17
\mathchardef\varthetag="7C23
\mathchardef\varpig="7C24
\mathchardef\varrhog="7C25
\mathchardef\varsigmag="7C26
\mathchardef\Omegag="7C0A
\mathchardef\Thetag="7C02
\mathchardef\Sigmag="7C06
\mathchardef\Deltag="7C01
\mathchardef\Phig="7C08
\mathchardef\Gammag="7C00
\mathchardef\Psig="7C09
\mathchardef\Lambdag="7C03
\mathchardef\Xig="7C04
\mathchardef\Pig="7C05
\mathchardef\Upsilong="7C07

\DeclareMathOperator*{\Spec}{Spec}

\def\ord{{\rm ord}}

\begin{document}

\title[Rational generating series]{On some
rational generating series occuring in arithmetic geometry}

\author{Jan Denef}
\address{University of Leuven, Department of Mathematics,
Celestijnenlaan 200B, 3001 Leu\-ven, Bel\-gium }
\email{ Jan.Denef@wis.kuleuven.ac.be}
\urladdr{http://www.wis.kuleuven.ac.be/wis/algebra/denef.html}

\author{Fran\c cois Loeser}

\address{{\'E}cole Normale Sup{\'e}rieure,
D{\'e}partement de math{\'e}matiques et applications,
45 rue d'Ulm,
75230 Paris Cedex 05, France
(UMR 8553 du CNRS)}
\email{Francois.Loeser@ens.fr}
\urladdr{http://www.dma.ens.fr/~loeser/}
\dedicatory{To Bernie's memory, with gratitude}

%\begin{abstract}We express the Lefschetz number of iterates of the
%monodromy of a function on a smooth complex algebraic variety 
%in terms of the Euler characteristic of a space of truncated arcs. 
%\end{abstract}

\maketitle

\section*{Introduction}

The main purpose 
of the present paper is to illustrate the
following motto:
``rational generating series occuring in arithmetic
geometry are motivic in nature''.
More precisely, consider a series
$F = \sum_{n \in \NN} a_n T^n$ 
with coefficients in $\ZZ$.
We shall say 
$F$ is motivic in nature
if there exists a series $F_{\rm mot} = \sum_{n \in \NN} A_n T^n$,
with coefficients $A_n$ in some Grothendieck ring of varieties,
or some Grothendieck ring of motives, 
such that $a_n$ is the number of rational points
of $A_n$ in some fixed finite field, for all $n \geq 0$.
Furthermore, we require $F_{\rm mot}$
to be canonically attached to $F$.
Of course, such a definition is somewhat
incomplete, since one can always take 
for $A_n$
the disjoint union of $a_n$ points.
In the present paper,
which is an update of a   talk by the second author
at the Conference
``Geometric Aspects of Dwork's theory''
that took place in
Bressanone in July 2001,
we consider the issue of being motivic
in nature for the following three types
of 
generating series: Hasse-Weil series, Igusa series and Serre series.
In section \ref{igc}, we consider the easiest case, that
of
Igusa type series, for which being motivic
in nature follows quite easily from 
Kontsevich's
theory of motivic integration as developed in \cite{jag} \cite{inv}.
The Serre case is more subtle. After a false start in
section \ref{sd1}, we explain in section \ref{sd2} how 
to deal with it by using
the work in \cite{def} on arithmetic motivic integration.
%The results in sections \ref{igc} to \ref{sd2},
%were all obtained in collaboration with Jan Denef.
Finally, in 
section \ref{dg}, we consider the case of
Hasse-Weil series, which still remains very much open.
Here there is a conjecture, which is due to
M. Kapranov \cite{kapranov} and can be traced back to
insights of Grothendieck
cf. p. 184 of \cite{gs}. 
Since a ``counterexample'' to the conjecture 
recently appeared \cite{ll}, we spend some time
to explain the dramatic  effects
of inverting the class of the affine line in
the Grothendieck group of varieties.
This gives us the opportunity of reviewing some interesting
recent work of Poonen \cite{poonen}, Bittner \cite{Bittner} and
Larsen and Lunts \cite{ll} and allows us
to  propose
a precised form of Kapranov's conjecture that escapes
Larsen and Lunts' counterexample.

\medskip It was one of Bernie's insights that most, if not all,
functions occuring in Number Theory should
be of geometric origin. So we hope the present contribution
will not be too inadequate as an homage to his memory.

\section{Conventions and preliminaries}\subsection{}In this paper, by a variety over a ring $R$, we mean
a reduced and separated scheme of finite type over $\Spec R$.

\subsection{}Let $A$ be a commutative ring. 
The ring of rational formal series with coefficients in $A$
is the smallest subring  of $A[[T]]$
containing $A [T]$ and stable under taking inverses (when they exist in
$A[[T]])$.

\section{Some classical generating series}

\subsection{}\label{A}Let $X$ be a variety over $\FF_q$. 
We set $N_n := \vert X (\FF_{q^n}) \vert$, for $n \geq 1$.

\begin{theorem}[Dwork \cite{dwork}]\label{thA}
The Hasse-Weil series
$$
Z (T)
:= \exp \Bigl ( \sum_{n \geq 1} \frac{N_n}{n} T^n \Bigr )
$$ 
is rational.
\end{theorem}

\subsection{}\label{B}Let $K$ be
a finite extension of $\QQ_p$ with ring of 
integers $\cO_K$ and uniformizing parameter
$\pi$.
Let $X$ be a variety over $\cO_K$ . We set
$\tilde N_n :=  \vert X (\cO_K / \pi^{n + 1}) \vert$, for $n \geq 0$.

\begin{theorem}[Igusa\cite{igusa}]\label{thB}
The series
$$
Q (T)
:=  \sum_{n \geq 0}  \tilde N_n T^n 
$$ 
is rational.
\end{theorem}

Strictly speaking
this result is due to Igusa \cite{igusa} in the hypersurface case and
to
Meuser in general \cite{meuser}. However, as mentioned in the
review MR 83g:12015
of \cite{meuser}, a trick by Serre allows to deduce
the general case from the hypersurface case.

\subsection{}\label{C}Let $K$ be
a finite extension of $\QQ_p$ with ring integers $\cO_K$ and uniformizing parameter
$\pi$. We keep the notations of \ref{B}.
For $n \geq 0$ we denote by $\overline N_n$ the cardinality of the image
of $X (\cO_K)$ in $X (\cO_K / \pi^{n + 1})$. In other words,
$\overline N_n$ is the number of points in
$X (\cO_K / \pi^{n + 1})$ (approximate solutions modulo $\pi^{n + 1}$)
that may be lifted to points in 
$X (\cO_K)$ (actual solutions in $\cO_K$). Clearly, $\overline N_n$
is finite. Furthermore, when $X$ is smooth, then $ \tilde N_n = \overline N_n$
for every $n$.

\begin{theorem}[Denef \cite{denef}]\label{thC}
The series
$$
P (T)
:=  \sum_{n \geq 0}  \overline N_n T^n 
$$ 
is rational.
\end{theorem}

\remark The problem of proving the analogue of Theorems \ref{thB}
and \ref{thC} when $K$ is a finite extension
of $\FF_q [[t]]$ still remains  very much
an open issue, but the level of difficulty
seems quite different for $Q (T)$ or $P (T)$.
While 
rationality of $Q (T)$ for
function fields would follow
using Igusa's proof once
Hironaka's strong form of resolution of singularities is known
in characteristic $p$, proving rationality 
of $P (T)$ for function fields would require completely new ideas,
since no general quantifier elimination 
Theorem is known, or even conjectured, in positive characteristic.

\section{Additive invariants of algebraic varieties}

\subsection{}\label{3.1}Let $R$ be a ring. We denote by ${\rm Var}_R$
the category of algebraic varieties over
$R$.
An additive invariant
$$
\lambda : 
{\rm Var}_R
\longrightarrow
S,
$$
with $S$ a ring, assigns to any $X$
in ${\rm Var}_R$ an element $\lambda (X)$ of $S$
such that
$$
\lambda (X) = \lambda (X')
$$
for $X \simeq X'$,
$$
\lambda (X) = \lambda (X') + \lambda (X \setminus X'),
$$
for $X'$ closed in $X$,
and
$$
\lambda (X \times X')
=
\lambda (X) \lambda (X')
$$
for every $X$ and $X'$.

Let us remark that additive invariants $\lambda$ naturally extend to take their values
on
constructible subsets  of algebraic varieties.

\subsection{Examples}
\subsubsection{Euler characteristic}Here $R = k$ is a field. When $k$ is a subfield
of $\CC$, the Euler characteristic
${\rm Eu} (X) := \sum_i (-1)^i {\rm rk} H^i_c (X (\CC), \CC)$
give rise to an additive invariant
${\rm Eu} : {\rm Var}_k \rightarrow \ZZ$. For general $k$, replacing
Betti cohomology with compact support
by $\ell$-adic cohomology with compact support,
$\ell \not= {\rm char} k$, one gets  an additive invariant
${\rm Eu}_{\ell} : {\rm Var}_k \rightarrow \ZZ$, which does not
depend on $\ell$.

\subsubsection{Hodge polynomial}Let us assume 
$R = k$ is a field of characteristic zero.
Then it follows from Deligne's Mixed Hodge Theory that
there is a unique 
additive invariant
$H : {\rm Var}_k \rightarrow \ZZ [u, v]$,
which assigns to a smooth projective variety $X$ over $k$
its usual Hodge polynomial
$$
H (u, v) := \sum_{p, q} (-1)^{p + q} h^{p,q} (X) u^p v^q ,
$$
with $h^{p,q} (X) = \dim H^q (X, \Omega^p_X)$ the $(p,q)$-Hodge number of $X$.

\subsubsection{Virtual motives}\label{vm}More generally, 
when $R = k$ is a field of characteristic zero, there exists by Gillet and Soul\'e \cite{GS},
Guillen and Navarro-Aznar \cite{GN},
a unique additive invariant $\chi_c : {\rm Var}_k \rightarrow K_0 ({\rm CHMot}_k)$,
which 
assigns to a smooth projective variety $X$ over $k$
the class of its Chow motive, where $ K_0 ({\rm CHMot}_k)$
denotes the Grothendieck ring of the category of Chow
motives over $k$ (with rational coefficients).

\subsubsection{Counting points} Counting points also
yields additive invariants.
Assume $k = \FF_q$, then
$N_n : X \mapsto \vert X (\FF_{q^n}) \vert$ gives rise to an additive invariant
$N_n : {\rm Var}_k \rightarrow \ZZ$.
Similarly, if $R$ is (essentially) of finite type over $\ZZ$, for every maximal
ideal $\gP$ of $R$ with finite residue field $k (\gP)$, we have 
an additive invariant
$N_{\gP} : {\rm Var}_R \rightarrow \ZZ$, which assigns to $X$
the cardinality of $(X \otimes k (\gP))(k (\gP))$.

\subsection{}There exists a universal
additive invariant
$[\_] : {\rm Var}_R \rightarrow K_0 ({\rm Var}_R)$ 
in the sense that 
composition with $[\_] $ gives a bijection between
ring morphisms $ K_0 ({\rm Var}_R) \rightarrow S$
and 
additive invariants ${\rm Var}_R \rightarrow S$.
The construction of $K_0 ({\rm Var}_R)$ is quite easy:
take the
free abelian group on isomorphism classes $[X]$
of objects of ${\rm Var}_R$ and mod out by the relation
$[X] = [X'] + [X \setminus X']$ for $X'$ closed in $X$.
The product is now defined by $[X] [X'] = [X \times X']$.

We shall denote by $\LL$ the class of the affine line $\AA^1_R$ in
$K_0 ({\rm Var}_R)$. An important role will be played by
the ring 
$\cM_R := K_0 ({\rm Var}_R)[\LL^{- 1}]$ obtained by localization with respect to the 
multiplicative set generated by $\LL$.
This construction is analogue to the construction
of the category of Chow motives from the category of effective
Chow motives by localization with respect to the Lefschetz motive.
(Remark that the morphism $\chi_c$
of \ref{vm} sends $\LL$ to the class of the Lefschetz motive.)

One should stress that very little is known about the structure
of the rings $K_0 ({\rm Var}_R)$ and $\cM_R$ even when $R$ is a field.
Let us just
quote a result by Poonen \cite{poonen}
saying that 
when $k$ is a field of characteristic zero the ring
$K_0 ({\rm Var}_k)$ is not a domain (we shall explain this result
with more details
in \S \ref{poo}). For instance, even for a field $k$, it is not known whether
the localization morphism
 $({\rm Var}_k) \rightarrow \cM_k$ is injective or not
(although the whole point of \S \ref{rc} relies on the guess it should not).

\remark In fact, the ring 
$K_0 ({\rm Var}_k)$
as well as the canonical morphism
$\chi_c : K_0 ({\rm Var}_k) \rightarrow
K_0 ({\rm CHMot}_k)$,
were already considered by Grothendieck in a letter to Serre dated
August 16, 1964, cf. p. 174 of \cite{gs}.

\section{Geometrization of $Q (T)$}\label{igc}

\subsection{}Let $k$ be a field. For every variety $X$ over $k$,
we denote by $\cL (X)$ the corresponding
space of arcs. It is a 
$k$-scheme, which satisfies
$$
\cL (X) (K) = X (K [[t]])
$$
for every field $K$ containg $k$.
More precisely $\cL (X)$ is defined
as the inverse limit
$\cL (X) := \limproj \cL_n (X)$,
where $ \cL_n (X)$ represents the functor from
$k$-algebras to sets sending a $k$-algebra $R$
to $X (R [[t]]/t^{n + 1} R [[t]])$.
We shall always consider $\cL (X)$ as endowed with its reduced structure.
We shall denote by $\pi_n$ the canonical morphism
$\cL (X) \rightarrow \cL_n (X)$.

\subsection{}We consider
the generating series
$$
Q_{\rm geom} (T) :=
\sum_{n \geq 0} \, [\cL_n (X)] \, T^n
$$
in $\cM_k [[T]]$.

\begin{theorem}[Denef-Loeser]\label{Qg}
Assume ${\rm char} k = 0$.
\begin{enumerate}
\item[1)]The series 
$Q_{\rm geom} (T)$
in $\cM_k [[T]]$ is rational of the form
$$
\frac{R (T)}{\prod (1 - \LL^a T^b)},
$$
with $R (T)$ in 
$\cM_k [T]$, $a$ in $\ZZ$ and $b$ in $\NN \setminus \{0\}$.
\item[2)] If $X$ is defined over
some number field $K$, then, for almost
all finite places $\gP$,
$$
N_{\gP} (Q_{\rm geom} (T)) = Q_{X \otimes \cO_{K_{\gP}}} (T).
$$
\end{enumerate}
\end{theorem}
Here we should explain what we mean by 
$N_{\gP} (Q_{\rm geom} (T)) $. For $X$ a variety over $K$,
$N_{\gP} (X)$ makes sense for almost all 
finite places $\gP$, by taking some model over $ \cO_{K_{\gP}}$.
Now we apply this termwise to the series $Q_{\rm geom} (T)$. 
This is possible since the series is rational by \ref{Qg} 1).
The Theorem  is proved for hypersurfaces in \cite{jag} , and the general case is similar and
may also be deduced from  general results in \cite{inv}
and \cite{def}.

\begin{proof}[Oversimplified sketch of proof the rationality]Let us first recall Igusa's proof of Theorem \ref{thB}
when
$X$ is an hypersurface defined by $f = 0$ in $\AA^m_{\cO_K}$.
The basic idea is to express the series
$Q (T)$ as the integral
$$
I (s) := \int_{\cO_K^m} \vert f \vert^s \vert dx \vert,
$$
up to trivial factors, with $T = q^{-s}$, $q$ the cardinality of the residue field.
Then one may use Hironaka's resolution of singularities to reduce
the computation of $I (s)$ to the case where $f = 0$ is locally given by monomials
for which direct 
calculation is easy.  

Our proof of the  rationality of $Q_{\rm geom} (T)$
follows similar lines. One express first our series as an integral, but
here
$p$-adic integration is replaced by \textbf{motivic integration}.
If $Y$ is a variety over $k$,
motivic integration assigns to certain subsets $A$ of the arc space $\cL (Y)$
a motivic measure $\mu (A)$ in $\cM_k$ (or sometimes, but
this will not be considered here,
a measure in a certain completion of $\cM_k$).
Then, to be able to use Hironaka's resolution of singularities to reduce
to the locally monomial case as  in Igusa's proof, we have to use the
fundamental
\textbf{change of variable formula} established in \S 3 of \cite{inv}.
\end{proof}

\section{Geometrization of $P(T)$: I}\label{sd1}

\subsection{}In view of the previous section, it is natural to consider
now the image $\pi_n (\cL (X))$ of $\cL (X)$ in $\cL_n (X)$.
Thanks to Greenberg's Theorem on solutions of polynomial systems in Henselian rings,
we know that $\pi_n (\cL(X))$ is a constructible subset$\cL_n (X)$, hence
we may consider its class $[\pi_n (\cL(X))]$ in $\cM_k$.
We consider the generating series
$$
P_{\rm geom} (T) :=
\sum_{n \geq 0} [\pi_n (\cL (X))] \, T^n
$$
in
$\cM_k [[T]]$.

\begin{theorem}[Denef-Loeser \cite{inv}]\label{Pg}
Assume ${\rm char} k = 0$.
The series 
$P_{\rm geom} (T)$
in $\cM_k [[T]]$ is rational of the form
$$
\frac{R (T)}{\prod (1 - \LL^a T^b)},
$$
with $R (T)$ in 
$\cM_k [T]$, $a$ in $\ZZ$ and $b$ in $\NN \setminus \{0\}$.
\end{theorem}

\begin{proof}[Oversimplified sketch of proof]Let us first recall
the strategy of the  proof \cite{denef}
of Theorem \ref{thC} in the $p$-adic case.
One reduces to the case where
$X$ is a closed subvariety of  $\AA^m_{\cO_K}$.
Then one  express the series
$P (T)$ as the integral
$$
J (s) := \int_{\cO_K^m} d (x, X)^s \vert dx \vert,
$$
up to trivial factors, with $T = q^{-s}$,
similarly as in Igusa's case, 
where 
$d (x, X)$ is the function ``distance to $X$''. Here an essential new 
feature
appears, the function
$d (x, X)$ being in general not
a polynomial function, but only a definable or semi-algebraic
function. Then one is  able to use
Macintyre's quantifier elimination Theorem \cite{angus}, a $p$-adic analogue
of Tarski-Seidenberg's theorem,
to prove rationality.

In the present setting our proof follows a similar pattern,
replacing $p$-adic integration by
motivic integration and the theory of 
$p$-adic semi-algebraic sets by 
a theory of 
$k [[t]]$-semi-algebraic sets
built off from a 
quantifier elimination Theorem due to Pas \cite{Pas}.
\end{proof}

\subsection{}When $X$ is defined over a number field
$K$, 
a quite natural guess would be, by
analogy with what we have seen
so far, that,
for almost
all finite places $\gP$,
$N_{\gP} (P_{\rm geom} (T)) = P_{X \otimes \cO_{K_{\gP}}} (T)$.
But such a statement cannot hold true.
This is due to the fact that,
in the very definition of $P (T)$, one is concerned in not considering
extensions of the residue field, while in the definition
of $P_{\rm geom} (T)$ extensions of the residue field $k$
are allowed.
To remedy this, one needs to be more careful about
rationality issues concerning the residue field, and for that purpose
it is convenient to introduce \textbf{definable subassignments}
as we do in the next section.

\section{Geometrization of $P(T)$: II}\label{sd2}

\subsection{Subassignments}
Fix a ring $R$. We denote
by $\textrm{Field}_R$
the category of $R$-algebras that are fields.
For an $R$-scheme $X$, we denote by $h_X$ the functor
which to a field $K$ in 
$\textrm{Field}_R$ assigns the set
$h_X (K) := X (K)$.
By a \textbf{subassignment} $h \subset h_X$
of 
$h_X$
we mean the datum, for every 
field $K$ in 
$\textrm{Field}_R$, of a subset $h (K)$ of 
$h_X (K)$. We stress that, contrarly to subfunctors,
no compatibility is required
between the various sets
$h (K)$.

All set theoretic constructions generalize in an obvious
way to the case of subassignments. For instance if $h$ and $h'$
are subassignments of $h_X$, then we denote by
$h \cap h'$ the subassignment $K \mapsto h (K) \cap h' (K)$, etc.

Also, if $\pi : X \rightarrow Y$ is a morphism
of $R$-schemes and $h$ is a 
subassignment of $h_X$, we define
the 
subassignment $\pi (h)$ of 
$h_Y$ by $\pi (h) (K) := \pi (h (K)) \subset h_Y (K)$.

\subsection{Definable subassignments}Let $R$
be a ring. 
By a ring formula $\varphi$  over $R$, we mean
a first order formula in the language of rings with coefficients
in $R$ and free variables $x_1, \dots x_n$.
In other words $\varphi$ is built out
from  boolean combinations
(``and'', ``or'', ``not'') of polynomial
equations over
$R$ and existential and universal quantifiers.

For example
$$
(\exists x)
(x^2 + x + y = 0 \;  \textrm{and} \; 4y \not= 1)
$$
is a ring 
formula over $\ZZ$ with free variable $y$.

To a  ring formula $\varphi$ over $R$ with
free variables $x_1, \dots x_n$
one assigns the subassignment $h_{\varphi}$
of $h_{\AA^n_R}$ defined by
\begin{equation}\label{hd}
h_{\varphi}( K) :=
\Bigl\{ (a_1, \dots, a_n) \in K^n\Bigm \vert
\varphi (a_1, \dots, a_n) \, \textrm{holds in } K \Bigr\} \subset
K^n = h_{\AA^n_R} (K).
\end{equation}

Such a subassignment 
of $h_{\AA^n_R}$ is called a \textbf{definable subassignment}.
More generally, using affine coverings, cf. \cite{def}, one defines 
definable subassignments of $h_X$ for $X$ a variety over $R$.

It is quite easy to show that if $\pi : X \rightarrow Y$
is an $R$-morphism of finite presentation,
$\pi (h)$ is a definable subassignment of $h_Y$
if $h$ is a definable subassignment of $h_X$.

In our situation, we are concerned with the
subassignment
$\pi (h_{\cL (X)}) \subset h_{\cL_n (X)}$.
Remark that
$\pi_n : \cL (X) \rightarrow \cL_n (X)$ is not of finite type.

Nevertheless, we have the following:

\begin{prop}[\cite{def}]
$\pi (h_{\cL (X)})$  is a definable subassignment of $h_{\cL_n (X)}$.
\end{prop}

\subsection{Formulas and motives}Let $k$ be a field of characteristic
zero.
It follows from \ref{vm} that we have a canonical morphism
$\chi_c :  K_0 ({\rm Var}_k) \rightarrow K_0 ({\rm CHMot}_k)$.
We shall denote by $K_0^{\rm mot} ({\rm Var}_k) $
the image of 
$K_0 ({\rm Var}_k)$
in
$K_0 ({\rm CHMot}_k)$ under this morphism.
Remark that the image of $\LL$
in
$K_0^{\rm mot} ({\rm Var}_k) $  is not a zero divisor since it is 
invertible in 
$K_0 ({\rm CHMot}_k)$.

Let us  explain now how to assign
in a canonical way to a ring formula 
$\varphi$ over $k$ an
element $\chi _c ([\varphi])$
of ${\rm K}_0^{\rm mot} ({\rm Var}_{k}) \otimes {\bf Q} $.

\subsection{}Let $\varphi$ be a formula over a number field $K$.
For almost all finite places $\gP$ with residue field
$k (\gP)$, one may extend the definition in
(\ref{hd}) to give a meaning
to $h_{\varphi}(k (\gP))$.
If $\varphi$ and $\varphi'$ are formulas over $K$,
we set $\varphi \equiv \varphi'$
if $h_{\varphi}(k (\gP)) =
h_{\varphi'}(k (\gP))$
for almost all finite places $\gP$.

It follows from a  fundamental result of 
J. Ax \cite{Ax} that
$\varphi \equiv \varphi'$
if and only if 
$h_{\varphi} (L) = h_{\varphi} (L')$
for every
pseudo-finite field $L$
containing $K$.
Let us recall that a
pseudo-finite field is an
infinite perfect field that has exactly one field extension of any given finite
degree,
and over which every geometrically irreducible variety has a rational point.
Historically, the above result of Ax was one of the main motivation
for introducing that notion.
One way of constructing
pseudo-finite fields
is by taking 
infinite ultraproducts of finite fields.

Let us now introduce
the Grothendieck
ring of formulas over $R$, $K_0 ({\rm Field}_R)$,
and
$K_0 ({\rm PFF}_R)$ the Grothendieck
ring of the 
theory of pseudo-finite fields over $R$.
The
ring $K_0 ({\rm Field}_R)$
(resp. $K_0 ({\rm PFF}_R)$)
is the group generated by symbols $[\varphi ]$, where $\varphi$
is any ring formula
over $R$, subject to the relations $[\varphi_1 \,{\rm or}\;
\varphi_2]=[\varphi_1]+[\varphi_2]- [\varphi_1 \, {\rm and}\;
\varphi_2]$, whenever $\varphi_1$ and $\varphi_2$ have the same free variables, and the
relations
$[\varphi_1] = [\varphi_2]$, whenever there exists a
ring formula $\psi$ over $k$ that, when interpreted in any 
field  (resp. any pseudo-finite field)
$K$ in ${\rm Field}_R$,
yields the graph of a bijection between the tuples of elements
of $K$ satisfying $\varphi_1$
and those satisfying $\varphi_2$. The ring multiplication 
is induced by the conjunction of formulas in disjoint sets of variables.

There is a canonical morphism
$$
K_0 ({\rm Field}_R) \longrightarrow K_0 ({\rm PFF}_R).
$$

We can now state
the following:

\begin{theorem}[Denef-Loeser \cite{def},\cite{icm}]\label{imp}Let $k$
be a field of characteristic zero.
There exists a unique ring morphism $$ \chi_c : {\rm K}_0
({\rm PFF}_k ) \longrightarrow {\rm K}_0^{\rm mot} ({\rm Var}_{k}  )
\otimes {\bf Q} $$
satisfying  the following two properties:
\begin{enumerate}\item[(i)] For any formula $\varphi$ which is a conjunction of polynomial equations
over $k$, the element $\chi_c ([\varphi])$ equals the
class in ${\rm K}_0^{\rm mot} ({\rm Var}_{k}  ) \otimes {\bf Q}$ of the
variety defined
by $\varphi$.
\item[(ii)] Let $X$ be a normal affine irreducible variety over $k$, $Y$ an unramified
Galois cover \footnote{Meaning that $Y$ is an integral \'etale scheme over
$X$ with $Y/G \cong X$, where $G$ is the group of all endomorphisms of $Y$ over $X$.}  of
$X$, and $C$ a cyclic subgroup of the Galois group G of $Y$ over $X$. For such data we
denote by
$\varphi_{Y,X,C}$ a ring formula, whose interpretation in any field $K$ containing
$k$, is the set of $K$-rational points on $X$ that lift to a geometric
point on $Y$ with decomposition group $C$ (i.e. the set of points on $X$ that lift to a $K$-rational
point of $Y/C$, but not to any $K$-rational point of $Y/C'$ with $C'$ a proper subgroup of
$C$). Then
\[
\chi_c ([\varphi_{Y,X,C} ]) = \frac{{\left| C \right|}}{{\left|
{{\rm N}_{G} (C)} \right|}}\chi_c ([\varphi_{Y,Y/C,C}
]),\]
where ${\rm N}_{G} (C)$ is the normalizer of $C$ in $G$.
\end{enumerate}
Moreover, when  $k$ is a number field,
for almost all finite places $\gP$,
$ N_{\gP} (\chi_c ([\varphi]))$ is
equal to the cardinality of $h_{\varphi} (k(\gP))$.
\end{theorem}

The above theorem is a variant of results in \S 3.4  of \cite{def}.
A sketch of proof is given in \cite{icm}.
\begin{proof}[Some ingredients in the proof]
Uniqueness uses quantifier elimination for
pseudo-finite fields (in terms of Galois stratifications,
cf. the work of
Fried and Sacerdote \cite{Fried Sacerdote}\cite[\S 26]{Fried Jarden}), 
from which 
it follows that ${\rm K}_0
({\rm PFF}_k )$ is generated as a group by
classes of
formulas of the form $\varphi _{Y,X,C}$.
Thus by (ii) we only have to
determine $\chi_c([\varphi _{Y,Y/C,C}])$, with $C$ a cyclic group.
But this follows
directly from the following recursion formula: 
\begin{equation}\label{rec}
\left| C \right|[Y/C] = \sum\limits_{A\,{\rm  subgroup\,  of }\, C} {\left| A \right|} \chi _c ([\varphi _{Y,Y/A,A}
]).
\end{equation}
This recursion formula is a direct consequence of (i), (ii), and the fact that
the formulas $\varphi _{Y,Y/C,A}$ yield a partition of $Y/C$.

The proof of existence
is based on work of del Ba{\~ n}o
Rollin and
Navarro Aznar
\cite{Rollin Aznar} who
associate to any representation over {\bf Q} of a finite group $G$ acting
freely on an affine variety $Y$ over $k$, an element in the Grothendieck group of Chow
motives over $k$. By linearity, we can hence associate to any {\bf Q}-central function
$\alpha$ on $G$ (i.e. a {\bf Q}-linear combination of characters of
representations of $G$ over {\bf Q}), an element $\chi_c(Y,\alpha )$ of that Grothendieck
group tensored with {\bf Q}. Using Emil Artin's Theorem, that any {\bf Q}-central function
$\alpha$ on $G$ is a {\bf Q}-linear combination of characters induced by trivial
representations of cyclic subgroups, one shows that
$\chi_c(Y,\alpha )\in {\rm K}_0^{\rm mot} ({\rm Var}_{k}  ) \otimes {\bf Q}$. For $X:=Y/G$
and $C$ any cyclic subgroup of $G$,
we define $\chi_c([\varphi _{Y,X,C}]) := \chi_c(Y,\theta)$, where $\theta$ sends $g\in G$ to 1
if the subgroup generated by $g$ is conjugate to $C$, and else to 0.
With some more work we prove that
the above definition of $\chi_c([\varphi _{Y,X,C}])$ extends by additivity to a
well-defined map
$ \chi _c : {\rm K}_0
({\rm PFF}_k ) \longrightarrow {\rm K}_0^{\rm mot} ({\rm Var}_{k}  ) \otimes {\bf Q} $. 
\end{proof}

Clearly $\chi_c (\varphi)$ depends only on $h_{\varphi}$ and the construction
easily extends by additivity
to 
definable subassignments of $h_X$, for
any variety $X$ over $k$. So,
to any such
definable subassignment $h$, we may associate $\chi_c (h)$ in
${\rm K}_0^{\rm mot} ({\rm Var}_{k}) \otimes \QQ$.

\begin{prop}[Denef-Loeser]Let $k$
be a field of characteristic zero.
For any 
definable subassignment $h$, ${\rm Eu} (\chi_c (h))$ belongs to $\ZZ$.
\end{prop}

\begin{proof}It is enough to show
that ${\rm Eu} (\chi_c (\varphi_{Y, X, C}))$ belongs to
$\ZZ$ for every $Y$, $X$ and $C$.
Consider first the case $C$
is the trivial subgroup $e$
of $G$. We have
$$\chi_c (\varphi_{Y, X, e}) = \frac{1}{|G|}
\chi_c (\varphi_{Y, Y, e}) =
\frac{1}{|G|}[Y].$$
It follows that 
$${\rm Eu} (\chi_c (\varphi_{Y, X, e})) =
\frac{1}{|G|}{\rm Eu} (Y)
= {\rm Eu} (X) \in \ZZ.$$
When 
$C$ is a non trivial cyclic subgroup of $G$, 
by induction on $|C|$, it follows
from
the recursion formula (\ref{rec})
that ${\rm Eu} (\chi_c (\varphi_{Y, X, C})) = 0$.
\end{proof}

\begin{example}Let $n$ be a integer $\geq 1$ and
assume $k$ contains
all $n$-roots of unity.
Consider the formula
$\varphi_n : (\exists y) (x = y^n \; \text{and} \;
x \not= 0)$ ; then $\chi_c (\varphi_n) = \frac{\LL -1}{n}$.
In particular
${\rm Eu} (\chi_c (\varphi_n)) = 0$
and
${H} (\chi_c (\varphi_n)) = \frac{uv -1}{n}$.
This example contradicts the example on page
430 line -2 of \cite{def}
(page 3 line 4 in the preprint) which is unfortunately incorrect.
\end{example}

\begin{remark}It is the place to correct
the following errors in the published version
of \cite{icm}.
On line 18 of
the third page, after the word ``motives'' one has to insert ``, and by killing
all {\bf L}-torsion''.
Once
this correction is made, it is easily
checked that
${\rm K}_0^{\rm mot} ({\rm Var}_{k})$ becomes the same
in the present paper and in \cite{icm}.
On line 6 of the eighth page, one has to delete the
last sentence.
\end{remark}

\subsection{The series $P_{\rm ar}$}

We now consider
the series
$$
P_{\rm ar} (T) := \sum_{n \geq 0} \chi_c (\pi_n (h_{\cL (X)})) \, T^n
$$
in
${\rm K}_0^{\rm mot} ({\rm Var}_{k}) \otimes \QQ$.

\begin{theorem}[Denef-Loeser \cite{def}]\label{Pa}
Assume ${\rm char} k = 0$.
\begin{enumerate}
\item[1)]The series 
$P_{\rm ar} (T)$
in ${\rm K}_0^{\rm mot} ({\rm Var}_{k})\otimes \QQ$ is rational of the form
$$
\frac{R (T)}{\prod (1 - \LL^a T^b)},
$$
with $R (T)$ in 
$({\rm K}_0^{\rm mot} ({\rm Var}_{k}) \otimes \QQ )[T]$, $a$ in $\ZZ$ and $b$ in $\NN \setminus \{0\}$.
\item[2)] If $X$ is defined over
some number field $K$, then, for almost
all finite places $\gP$,
$$
N_{\gP} (P_{\rm ar} (T)) = P_{X \otimes \cO_{K_{\gP}}} (T).
$$
\end{enumerate}
\end{theorem}

%\begin{remark}In fact one can show that
%$P_{\rm ar} (T)$ is already rational in the ring
%$K_0 ({\rm Form}_k) [[T]]$.
%\end{remark}

In the proof of Theorem \ref{Pa}, one uses in an essential way
\textbf{arithmetic motivic integration} 
a variant of motivic integration developed in \cite{def}.
The specialization statement 
2) in Theorem \ref{Pa} is a special case of the following
results, which states that
``natural $p$-adic integrals are motivic''.

\begin{theorem}[Denef-Loeser \cite{def}]\label{int}
Let $K$ be a number field. Let $\varphi$ be
a first order formula in the language
of \textbf{valued} rings with coefficients
in $K$ and free variables 
$x_1$, \dots, $x_n$. Let $f$ be a polynomial
in $K [x_1, \dots, x_n]$.
For 
$\got P$ a finite place of $K$, denote by 
$K_{\got P}$ the completion of $K$
at $\got P$. Then there exists a \textbf{canonical} motivic integral
which specializes to
$$
\int_{h_{\varphi} (K_{\got P})} |f|_{\got P}^s |d x|_{\got P}
$$
for almost all finite places $\got P$.
\end{theorem}

The formulation here is somewhat unprecise and we refer to 
\cite{def} for details.
Let us just say that formulas in the language
of valued rings include 
expressions like
$ \ord a \leq \ord b$ or $\ord c \equiv b \mod e$.

\subsection{}Theorem \ref{int} applies in particular to integrals occuring
in $p$-adic harmonic analysis, like orbital integrals.
This has led recently Tom Hales to speculate that
many of the basic objects in
representation theory are motivic in nature and to
develop 
a beautiful conjectural program aiming to the determination
of the (still conjectural) virtual Chow motives
that control the behavior of orbital integrals
and leading to a motivic fundamental lemma \cite{hales}\cite{gh}.

\section{Geometrization of $Z(T)$}\label{dg}

\subsection{}Let $k$ be a field and
let $X$ be a variety over 
$k$. For $n \geq 0$, we denote bt $X^{(n)}$ the $n$-fold symmetric
product of $X$, i.e. the quotient of the cartesian product
$X^n$ by the symmetric group of $n$ elements. Note that
$X^{(0)}$ is isomorphic to $\Spec k$.

Following Kapranov \cite{kapranov}, we define the motivic
zeta function of $X$ as the power series
$$
Z_{\rm mot} (T) :=  \sum_{n = 0}^{\infty}
[X^{(n)}] \, T^n
$$
in $K_0 (\Var_k) [[T]]$.

Also, when $\alpha : K_0 (\Var_k) \rightarrow A$
is a morphism of rings, we denote by
$Z_{\rm mot, \alpha} (T) $ the
power series $\sum_{n = 0}^{\infty}
\alpha([X^{(n)}])  \, T^n$
in $A [[T]]$.
We shall write $\LL$ for $\alpha (\LL)$.

\begin{prop}If $k = \FF_q$, and we write
$N (S) = N_1 (S) = \vert S (k) \vert$ for $S$ a variety over $k$,
then 
$Z_{\rm mot, N} (T) $
is equal to the Hasse-Weil zeta function considered in \ref{A}.
\end{prop}

\begin{proof}Rational points of $X^{(n)}$ over $k$
correspond to degree $n$
effective zero cycles of $X$,
hence the result follows from the usual inversion formula
between the number of 
effective zero cycles of given degree on $X$ 
and the number of rational points of $X$ over finite extensions of $k$.
\end{proof}

In his paper \cite{kapranov}, Kapranov
proves the following result:

\begin{theorem}[Kapranov \cite{kapranov}]Let $X$ be a smooth projective irreducible curve of genus
$g$. Let 
$\alpha : K_0 (\Var_k) \rightarrow A$
be a morphism of rings with 
$A$ a field, such that $\LL$ is non zero in $A$.
Assume also there exists a degree 1 line bundle on $X$.
Then:
\begin{enumerate}
\item[1)]The series
$Z_{\rm mot, \alpha}  (T)$ is rational. It is the quotient of a polynomial of degree $2g$
by $(1 -T) (1 -\LL T)$.
\item[2)]The function
$Z_{\rm mot, \alpha}  (T)$ satisfies the functional equation
$$
Z_{\rm mot, \alpha}  (\LL^{-1}T^{-1})
= \LL^{1 -g} T^{2- 2g}
Z_{\rm mot, \alpha}  (T).
$$
\end{enumerate}
\end{theorem}

The proof follows the lines of F. K. Schmidt's  classical proof \cite{fk}
of rationality and functional equation
for the Hasse-Weil zeta function of a smooth projective curve.

In the same paper, Kapranov states `` it is natural to expect
that the motivic zeta functions are rational and satisfy similar functional
equations?".

\remark A generating series similar to $Z_{\rm geom}$ and the question of
its rationality
were already considered by Grothendieck in a letter to Serre dated
September 24, 1964, cf. p. 184 of \cite{gs}.

\subsection{Stable birational invariants}
We now give a new presentation by generators and relations of 
$K_0 (\Var_k)$
due to F. Bittner \cite{Bittner}.
We denote
by  $K_0^{\rm bl} (\Var_k)$
the quotient of the free abelian group
on isomorphism classes of smooth proper 
varieties over $k$
by the relation
$$
[{\rm Bl}_Y X] - [E] =  [X] - [Y],
$$
for $Y$ and $X$ smooth proper irreducible over $k$,
$Y$ closed in $X$, ${\rm Bl}_Y X$ the blowup of $X$ with center
$Y$ and $E$ the exceptional divisor in ${\rm Bl}_Y X$.
As for $K_0 (\Var_k)$,
cartesian product induces a product on 
$K_0^{\rm bl} (\Var_k)$ which endowes it with a ring structure.
There is a canonical ring morphism
$
K_0^{\rm bl} (\Var_k) \rightarrow K_0 (\Var_k),
$
which sends $[X]$ to $[X]$.

\begin{theorem}[Bittner \cite{Bittner}]\label{bi}Assume $k$
is of characteristic zero. The
canonical ring morphism
$$
K_0^{\rm bl} (\Var_k) \rightarrow K_0 (\Var_k)
$$
is an isomorphism.
\end{theorem}

The proof is based on Hironaka resolution of singularities
and the weak factorization Theorem
of Abramovich, Karu, Matsuki and
W{\l}odarczyk \cite{wf}.

One deduces easily the following result, first proved by
Larsen and Lunts \cite{ll}.

\begin{cor}[Larsen and Lunts \cite{ll}]\label{lalu}
Let us assume $k$ is algebraically closed
of characteristic zero.
Let $A$ be the monoid of isomorphism classes
of smooth proper irreducible
varieties over $k$ and let
$\Psi : A \rightarrow G$ be a morphism
of commutative monoids such that
\begin{enumerate}
\item [1)]If $X$ and $Y$ are
birationally equivalent
smooth proper irreducible
varieties over $k$, then $\Psi ([X]) = \Psi ([Y])$.
\item [2)]$\Psi ([\PP^n_k]) = 1$.
\end{enumerate}
Then there exists a unique morphism a rings
$$
\Phi : K_0 (\Var_k) : \longrightarrow \ZZ [G]
$$
such that $\Phi ([X]) = \Psi ([X])$
when $X$ is smooth proper irreducible.
\end{cor}

We assume from now on that $k$ is algebraically closed
of characteristic zero.
We denote by ${\rm SB}$ the monoid
of 
equivalence classes
of smooth proper irreducible
varieties over $k$ under stably birational
equivalence \footnote{$X$ and $Y$ are called
stably birational if $X \times \PP^r_k$ is
birational to $Y \times \PP^s_k$ for some $r, s \geq 0$.}. It follows
from
Corollary \ref{lalu} that there exists a  
universal stable birational
invariant 
$$
\Phi_{\rm SB} : K_0 (\Var_k) : \longrightarrow \ZZ [{\rm SB}].
$$

\begin{prop}[Larsen and Lunts \cite{ll}]\label{lalu2}The kernel of the morphism
$\Phi_{\rm SB} : K_0 (\Var_k) : \rightarrow \ZZ [{\rm SB}]$
is the principal
ideal generated by $\LL = [\AA^1_k]$.
\end{prop}

\begin{proof}[Sketch of proof]It is clear that $\LL$ lies in the kernel
of 
$\Phi_{\rm SB}$. Conversely,
take $\alpha = \sum_{1 \leq i \leq r} [X_i] -
\sum_{1 \leq j \leq s} [Y_j]$ in 
the kernel
of 
$\Phi_{\rm SB}$,
with $X_i$ and $Y_j$ smooth, proper and irreducible.
Since 
$\sum_{1 \leq i \leq r} [X_i] = 
\sum_{1 \leq j \leq s} [Y_j]$ in $\ZZ [\rm SB]$, $r = s$ and, after renumbering
the $X_i$'s, we may assume 
$X_i$ is stably birational to $Y_i$ for every $i$. 
Hence it is enough to show that if $X$ and $Y$ are
smooth, proper and irreducible stably birationally equivalent,
then $[X] - [Y]$ belongs
to $\LL K_0 (\Var_k)$. Since $[X] - [\PP^r_k \times X]$
belongs
to $\LL K_0 (\Var_k)$, we can even assume $X$ and $Y$ are birationally
equivalent and then the result follows easily from the
weak factorization Theorem.
\end{proof}

\subsection{Back to Poonen's result}\label{poo}As promised, we shall now give some
explanations concerning the proof of Poonen's Theorem \ref{pt}.

\begin{klem}[Poonen \cite{poonen}]\label{kl}Let $k$ be a field of characteristic
zero.
There exists abelian varieties $A$ and $B$ over $k$ such that
$A \times A$ is isomorphic to $B \times B$ but
$A_{\overline k} \not\cong B_{\overline k}$.
\end{klem}

The proof relies on the following lemma:

\begin{lem}[Poonen \cite{poonen}]Let $k$ be a field of characteristic
zero. There exists an 
abelian variety $A$ over $k$ 
such that ${\rm End}_k (A) = {\rm End}_{\overline k} (A)
\cong \cO$, with $\cO$ the ring of integers of a number field of class
number 2.
\end{lem}

When $k = \CC$, one may take $A$ an elliptic curve with complex multiplication
by $\ZZ [\sqrt{-5}]$. The general case is much more involved and 
necessitates the use of
modular forms and Eichler-Shimura Theory as well as some table checking,
see \cite{poonen}.

Let us now explain how Poonen deduces from the Key-Lemma the following:

\begin{theorem}\label{pt}The ring $K_0 (\Var_k)$ is not a domain, for 
$k$ a field of characteristic
zero.
\end{theorem}

\begin{proof}Take $A$ and $B$ as in the
Key-Lemma. We have
$([A] + [B]) ([A] - [B]) = 0$ in 
$K_0 (\Var_k)$.
To check that 
$[A] + [B]$ and  $[A] - [B]$ are nonzero in $K_0 (\Var_k)$,
it is
enough to check that they have a nonzero image
under the composition
$$
K_0 (\Var_k) \rightarrow 
K_0 (\Var_{\overline k}) \rightarrow 
\ZZ [{\rm SB}_{\overline k}] \rightarrow 
\ZZ [{\rm AV}_{\overline k}],
$$
where ${\rm AV}_{\overline k}$ is the monoid of
isomorphism classes of abelian varities over $\overline k$
and the last morphism is induced by the Albanese functor
assigning to a smooth irreducible variety its Albanese variety
(which is indeed a stable birational invariant).
To conclude we just have to remark that the 
Albanese variety of an abelian variety is equal to itself.
\end{proof}

\begin{remark}Poonen's proof does not tell us anything about
zero divisors in $\cM_k$. Indeed, it relies on the use of stable
birational invariants, and after  inverting $\LL$ no (non trivial)
such  invariant is left.
\end{remark}

\subsection{Non rationality results}
Larsen and Lunts
proved \cite{ll} the following 
non rationality  Theorem:

\begin{theorem}[Larsen-Lunts]\label{nr}
Let $X$ be a smooth proper complex irreducible
surface with geometric genus
$p_g (X) \geq 2$.
Then there exists  a
morphism
$\alpha :  K_0 (\Var_k) \rightarrow  F$, with $F$
a field, such that the
zeta function $Z_{\rm mot, \alpha}$
attached to $X$ is not rational.
\end{theorem}

\begin{proof}[Some ideas from the proof]Larsen and Lunts consider,
for $X$ a smooth proper complex irreducible variety of
dimension $d$,
the polynomial
$\Psi_h (X) := \sum_{1 \leq i \leq d} h^{i, 0} t^i$.
Remark $p_g (X)$ is the leading coefficient of
$\Psi_h (X)$.
It is well known $\Psi_h$ is a stable birational invariant,
hence by Corollary \ref{lalu} it gives rise
to a ring morphism
$$
\Phi_h : K_0 (\Var_k) \longrightarrow
\ZZ [C],
$$
with $C$
the multiplicative 
monoid of polynomials in $\ZZ [t]$ with positive leading coefficient.
Larsen and Lunts  show that  $\ZZ [C]$
is a domain and take for $\alpha$ the composition of
$\Phi_h$ with the localization morphism
from  $\ZZ [C]$ to its fraction field
$F$.
A key ingredient
in the proof is the fact
that,  for  smooth projective surfaces
with geometrical genus $r \geq 2$, 
$p_g (X^{(n)}) ={\binom{r + n - 1}{r- 1}}$.
Here 
$p_g (X^{(n)})$ is by definition
the geometric genus
of any smooth proper variety birational to 
$X^{(n)}$. The Hilbert scheme
$X^{[n]}$ parametrizing closed zero-dimensional subschemes
of length $n$ of $X$ is such a variety and 
it follows from results of
G\"ottsche and Soergel \cite{goso}
that
$p_g (X^{[n]}) ={\binom{r + n - 1}{r- 1}}$.
Some
more work is needed in order to deduce
non rationality for the series $Z_{\rm mot, \alpha}$.
\end{proof}

\begin{remark}Since the first version of the present
paper was written, 
Larsen an Lunts wrote 
a
very interesting sequel \cite{ll2}
of \cite{ll}.
\end{remark}

\subsection{Rationality conjectures}\label{rc}In view of Theorem \ref{nr},
one cannot hope for the series $Z_{\rm mot}$ to be in general rational
in $K_0 (\Var_k) [[T]]$.

Nevertheless, one can still believe in 

\begin{sconjecture}\label{sc}Let $X$ be a variety over a field $k$.
Then the series $Z_{\rm mot}$ attached to $X$ is
rational in $\cM_k [[T]]$.
\end{sconjecture}

\begin{wconjecture}Let $X$ be a variety over a field $k$.
Then, for every
morphism $\alpha :  \cM_k \rightarrow F$, with $F$ a field,
the series $Z_{\rm mot, \alpha }$ attached to $X$ is
rational in $F [[T]]$.
\end{wconjecture}

\begin{remarks}
\begin{enumerate}
\item[1)]A posteriori it is not so surprising that we have to 
invert $\LL$ in order for rationality to conjecturally hold. Indeed,
the guess that the motivic series should be rational comes from
analogy with Dwork's Theorem \ref{thA}. But counting rational points is certainly not
a birational invariant!
\item[2)]When $X$ is smooth and proper, one can conjecture strong and weak forms
of functional equations.
%\item[3)]In case $\cM_k$ happened to be a domain, the strong and weak forms
%of the conjecture would be equivalent by Lemma \ref{tr}.
%\item[4)]In case the canonical morphism
%$K_0 (\Var_k) \rightarrow \cM_k$ would be injective (which I guess not to happ%en),
%Conjecture \ref{sc}
%could not be true,
%due to Lemma \ref{tr} and
%Theorem \ref{nr}.
\end{enumerate}
\end{remarks}

\bibliographystyle{amsplain}

\end{document}